\documentclass[12pt,a4paper]{article}

\usepackage[usenames]{color}
\usepackage{amssymb}
\usepackage{amsmath}
\usepackage{amsthm}
\usepackage{amsfonts}
\usepackage{amscd}
\usepackage{graphicx}

\usepackage[colorlinks=true,
linkcolor=webgreen,
filecolor=webbrown,
citecolor=webgreen]{hyperref}

\definecolor{webgreen}{rgb}{0,.5,0}
\definecolor{webbrown}{rgb}{.6,0,0}
\usepackage{color}
\usepackage{fullpage}
\usepackage{float}

\usepackage{cases}

\usepackage{graphics,amsmath,amssymb}
\usepackage{amsthm}
\usepackage{amsfonts}
\usepackage{latexsym}
\usepackage{enumitem}

\begin{document}

\theoremstyle{plain}
\newtheorem{theorem}{Theorem}
\newtheorem{remark}{Remark}
\newtheorem{proposition}{Proposition}
\newtheorem{example}[remark]{Example}
\newtheorem{lemma}{Lemma}
\newtheorem{definition}{Definition}
\newtheorem{corollary}{Corollary}

\begin{center}
\vskip 1cm{\Large\bf A remark on an explicit formula for
\vskip .06in the sums of powers of integers}
\vskip .3in \large Jos\'{e} L. Cereceda \\
{\normalsize Collado Villalba, 28400 (Madrid), Spain} \\
\href{mailto:jl.cereceda@movistar.es}{\normalsize{\tt jl.cereceda@movistar.es}}
\end{center}

\begin{abstract}
Recently, E. Samsonadze (arXiv:2411.11859v1) has given an explicit formula for the sums of powers of integers $S_k(n) = 1^k +2^k +\cdots + n^k$. In this short note, we show that Samsonadze's formula corresponds to a well-known formula for $S_k(n)$ involving the Stirling numbers of the second kind.
\end{abstract}

\vspace{4mm}

For integers $k \geq 0$ and $n \geq 1$, let $S_k(n)$ denote the power sum $S_k(n) = 1^k + 2^k + \cdots + n^k$. In a recent paper, Samsonadze \cite[Theorem 2.4]{samson} derived the following explicit formula for $S_k(n)$ (in our notation):
\begin{equation}\label{samson}
S_k(n) = (-1)^k \sum_{j=0}^k a_{k,j} n(n+1) \ldots (n+j),
\end{equation}
where
\begin{equation}\label{samson2}
a_{k,j} = \frac{1}{j+1} \sum_{i=0}^j \frac{(-1)^i i^k}{i! (j-i)!}.
\end{equation}
(Note that we are writing $a_{k,j}$ instead of merely $a_j$ because $a_{k,j}$ depends explicitly on both $k$ and $j$. Note also that we make the index $j$ start at $j=0$.) By combining \eqref{samson} and \eqref{samson2}, the above formula for $S_k(n)$ can be written in terms of binomial coefficients as follows:
\begin{equation}\label{samson3}
S_k(n) = (-1)^k \sum_{j=0}^k j! \binom{n+j}{j+1} \frac{1}{j!} \sum_{i=0}^j (-1)^i \binom{j}{i} i^k.
\end{equation}

On the other hand, as is well known, the Stirling numbers of the second kind $S(k,j)$ have the explicit formula (see, e.g., \cite[Equation (4)]{boya})
\begin{equation*}
S(k,j) = \frac{1}{j!} \sum_{i=0}^j  (-1)^{j-i} \binom{j}{i} i^k.
\end{equation*}
Therefore, formula \eqref{samson3} (or formula \eqref{samson} by the way) can in turn be written as
\begin{equation}\label{yo}
S_k(n) = \sum_{j=0}^k (-1)^{k-j} j! S(k,j) \binom{n+j}{j+1},
\end{equation}
for all $k \geq 0$. Furthermore, the coefficients $a_{k,j}$ in \eqref{samson2} can be expressed alternatively as
\begin{equation*}
a_{k,j} = \frac{(-1)^j}{j+1} S(k,j), \quad j =0,1,\ldots, k.
\end{equation*}
In particular, $a_{k,k} = \frac{(-1)^k}{k+1}$ for all $k \geq 0$ and $a_{k,1} = -\frac{1}{2}$ for all $k \geq 1$.

Formula \eqref{yo} is certainly a well-known formula for $S_k(n)$ involving the Stirling numbers of the second kind; see, e.g., Equation (7.5) of \cite{gould}, Equation (21) of \cite{howard}, the second formula at the top of p.\ 285 of \cite{knuth}, and Equation (15.27) of \cite{quain}. Some more information about formula \eqref{yo} can be found in section 4 of \cite{cere}.

Moreover, it is to be mentioned that formula \eqref{yo} is closely related to the companion formula
\begin{equation}\label{yo2}
S_k(n) = \sum_{j=0}^k j! S(k,j) \binom{n+1}{j+1}, \quad k \geq 1.
\end{equation}
In fact, formula \eqref{yo} can be obtained from \eqref{yo2} (and vice versa) by making the transformation $n \to -n-1$ and using the symmetry property of the power sum polynomial $S_k(n)$. Several applications of this transformation are given in \cite{cere2}.

Finally, it is readily verified that formula \eqref{yo} can be factorized as
\begin{equation*}
S_k(n) = 2S_1(n) \sum_{j=1}^k (-1)^{k-j} \frac{(j-1)!}{j+1} S(k,j) \binom{n+j}{j-1},
\quad k \geq 1,
\end{equation*}
which is equivalent to Theorem 2.5 of \cite{samson}.


\vspace{-1mm}

\begin{thebibliography}{9}

\small{
\bibitem{boya} K. N. Boyadzhiev, Close encounters with the Stirling numbers of the second kind, Mathematics Magazine 85(4) (2012), pp. 252--266.
\vspace{-2mm}

\bibitem{cere} J. L. Cereceda, Polynomial interpolation and sums of powers of integers, International Journal of Mathematical Education in Science and Technology 48(2) (2017), pp. 267--277.
\vspace{-2mm}

\bibitem{cere2} J. L. Cereceda, Dual recursive formulas for the sums of powers of integers, Far East Journal of Mathematical Education 26(2) (2024), pp. 111-121.
\vspace{-2mm}

\bibitem{gould} H. W. Gould, Evaluation of sums of convolved powers using Stirling and Eulerian numbers, The Fibonacci Quarterly 16(6) (1978), [Part 1], pp. 488--497.
\vspace{-2mm}

\bibitem{howard} F. T. Howard, Sums of powers of integers, Mathematical Spectrum 26(4) (1993/1994), pp. 103--109.
\vspace{-2mm}

\bibitem{knuth} D. E. Knuth, Johann Faulhaber and sums of powers, Mathematics of Computation 61(203) (1993), pp. 277--294.
\vspace{-2mm}

\bibitem{quain} J. Quaintance and H. W. Gould, Combinatorial Identities for Stirling Numbers. The Unpublished
Notes of H W Gould. World Scientific Publishing, Singapore, 2016.
\vspace{-2mm}

\bibitem{samson} E. Samsonadze, On sums of powers of natural numbers. Preprint (2024), available at \url{https://arxiv.org/abs/2411.11859v1}}

\end{thebibliography}
\end{document}